\documentclass[a4paper]{article}     
\usepackage{theorem}                 
\usepackage{amssymb}
\usepackage{epsfig}
\newtheorem{theorem}{Theorem}

\newtheorem{proposition}{Proposition}

\newcommand{\CC}{{\mathbb{C}}}
\newcommand{\HH}{{\mathbb{H}}}

\newcommand{\RR}{{\mathbb{R}}}
\newcommand{\ZZ}{{\mathbb{Z}}}

\newcommand{\NN}{{\mathbb{N}}}

\begin{document}

\title{Strange duality, mirror symmetry, and the Leech lattice}
\author{Wolfgang Ebeling}
\date{}

\maketitle

\rightline{\large{\it Dedicated to Terry Wall}}

\begin{abstract}
We give a survey on old and new results concerning Arnold's strange duality. We
show that most of the features of this duality continue to hold
for the extension of it discovered by C.~T.~C.~Wall and the author. The
results include relations to mirror symmetry and the Leech lattice.
\end{abstract}

\section*{Introduction}

More than 20 years ago, V.~I.~Arnold \cite{Arnold75} discovered a strange
duality
among the 14 exceptional unimodal hypersurface singularities. A beautiful
interpretation of this duality was given by H.~Pinkham \cite{Pinkham77} and
independently by I.~V.~Dolgachev and V.~V.~Nikulin \cite{DN77, Dolgachev82}.
I.~Nakamura related this duality to the Hirzebruch-Zagier duality of cusp
singularities
\cite{Nakamura80, Nakamura81}.

In independent work in early 1982, C.~T.~C.~Wall and the author discovered an
extension of this duality embracing on one hand series of bimodal singularities
and on the other hand, complete intersection surface singularities in $\CC^4$
\cite{EW85}. We showed that this duality also corresponds to Hirzebruch-Zagier
duality of cusp singularities.

Recent work has aroused new interest in Arnold's strange duality. It was
observed by several authors (see
\cite{Dolgachev95} and the references there) that Pinkham's interpretation of
Arnold's original strange duality can be considered as part of a two-dimensional
analogue of the mirror symmetry of families of Calabi-Yau threefolds. Two years
ago, K.~Saito
\cite{Saito94} discovered a new feature of Arnold's strange duality involving
the characteristic polynomials of the monodromy operators of the singularities
and he found a connection with the characteristic polynomials of automorphisms
of the famous Leech lattice. Only shortly after, M.~Kobayashi
\cite{Kobayashi95} found a duality of the weight systems associated to the 14
exceptional unimodal singularities which corresponds to Arnold's strange
duality. He also related it to mirror symmetry.

In this paper we first review these results.

Then we consider our extension of this duality and examine which of the
newly discovered features continue to hold. It turns out that with a suitable
construction, Pinkham's interpretation can be extended to a larger class of
singularities. In this way, one obtains many new examples of mirror symmetric
families of K3 surfaces.  We also associate characteristic polynomials to the
singularities involved in our extension of the duality and show that Saito's
duality continues to hold. Moreover, in this way we can realize further
characteristic polynomials of automorphisms of the Leech lattice. The connection
with the Leech lattice seems to be rather mysterious. We discuss some facts
which
might help to understand this connection. We conclude with some open questions.

We thank the referee for his useful comments.

\section{Arnold's strange duality}

We first discuss Arnold's original strange duality among
the 14 exceptional unimodal hypersurface singularities.

We recall Dolgachev's construction \cite{Dolgachev74, Dolgachev75} (see also
\cite{Looijenga83}) of these singularities. Let $b_1 \leq b_2 \leq b_3$ be
positive integers such that
$\frac{1}{b_1} + \frac{1}{b_2} + \frac {1}{b_3} < 1$. Consider the upper half
plane $\HH = \{ x+iy \in \CC | y >0\}$ with the hyperbolic metric
$\frac{1}{y^2}(dx^2+dy^2)$ and a solid triangle
$\Delta
\subset \HH$ with angles $\frac{\pi}{b_1}$, $\frac{\pi}{b_2}$,
$\frac{\pi}{b_3}$. Let $\Sigma$ be the subgroup of the group of isometries of
$\HH$ generated by the reflections in the edges of $\Delta$, and let $\Sigma_+$
be the subgroup of index 2  of orientation preserving isometries. Then
$\Sigma_+ \subset {\rm PSL}_2(\RR)$ and $\Sigma_+$ acts linearly on the  total
space $T\HH$ of the tangent bundle on $\HH$. The inclusion $\HH \subset T\HH$
as zero section determines an inclusion $\HH / \Sigma_+ \subset T\HH /
\Sigma_+$ of orbit spaces. Collapsing $\HH / \Sigma_+$ to a point yields a
normal surface singularity $(X,x_0)$. This singularity is called a {\em
triangle singularity}. The numbers $b_1$, $b_2$, $b_3$ are called the {\em
Dolgachev numbers} ${\rm Dol}(X)$ of the singularity.
The scalar multiplication in the fibres of the tangent bundle $T\HH$ induces a
good $\CC^\ast$-action on $X$. A resolution of the singularity $(X,x_0)$ can be
obtained by the methods of \cite{OW71}. A minimal good resolution consists of a
rational curve of self-intersection number $-1$ and three rational curves of
self-intersection numbers $-b_1$, $-b_2$, and  $-b_3$ respectively intersecting
the exceptional curve transversely.

By \cite{Dolgachev74}, for exactly 14 triples
$(b_1, b_2, b_3)$ the singularity $(X,x_0)$ is a hypersurface singularity.
Thus it can be given by a function germ
$f: (\CC^3,0) \to (\CC,0)$ where  $f$ is weighted
homogeneous with weights $w_1$, $w_2$, $w_3$ and degree $N$. The corresponding
weighted homogeneous functions, weights and degrees are indicated in
Table~\ref{Table1}. It turns out that these singularities are unimodal and one
gets in this way exactly the 14 exceptional unimodal hypersurface singularities
in Arnold's classification \cite{Arnold75}. (The equations in
Table~\ref{Table1} are obtained by setting the module equal to zero.)

\begin{table}\centering
\caption{The 14 exceptional unimodal singularities} \label{Table1}
\begin{tabular}{|c|c|c|c|c|c|c|c|c|}  \hline
Name   & Equation & $N$ & Weights & Dol & Gab & $\mu$ & $d$  &
Dual
\\ \hline
$E_{12}$ & $x^7+y^3+z^2$ & 42 & 6 14 21 & 2 3 7 & 2 3 7 & 12 & 1
 & $E_{12}$ \\ \hline
$E_{13}$ & $x^5y+y^3+z^2$ & 30 & 4 10 15 & 2 4 5 & 2 3 8 & 13 & $-2$
 & $Z_{11}$ \\ \hline
$E_{14}$ & $x^8+y^3+z^2$ & 24 & 3 8 12 & 3 3 4 & 2 3 9 & 14 & 3
 & $Q_{10}$ \\ \hline
$Z_{11}$ & $x^5+xy^3+z^2$ & 30 & 6 8 15 & 2 3 8 & 2 4 5 & 11 & $-2$
 & $E_{13}$ \\ \hline
$Z_{12}$ & $x^4y+xy^3+z^2$ & 22 & 4 6 11 & 2 4 6 & 2 4 6 & 12 & 4
 & $Z_{12}$ \\ \hline
$Z_{13}$ & $x^6+xy^3+z^2$ & 18 & 3 5 9 & 3 3 5 & 2 4 7 & 13 & $-6$
 & $Q_{11}$ \\ \hline
$Q_{10}$ & $x^4+y^3+xz^2$ & 24 & 6 8 9 & 2 3 9 & 3 3 4 & 10 & 3
 & $E_{14}$ \\ \hline
$Q_{11}$ & $x^3y+y^3+xz^2$ & 18 & 4 6 7 & 2 4 7 & 3 3 5 & 11 & $-6$
 & $Z_{13}$ \\ \hline
$Q_{12}$ & $x^5+y^3+xz^2$ & 15 & 3 5 6 & 3 3 6 & 3 3 6 & 12 & 9
 & $Q_{12}$ \\ \hline
$W_{12}$ & $x^5+y^4+z^2$ & 20 & 4 5 10 & 2 5 5 & 2 5 5 & 12 & 5
 & $W_{12}$ \\ \hline
$W_{13}$ & $x^4y+y^4+z^2$ & 16 & 3 4 8 & 3 4 4 & 2 5 6 & 13 & $-8$
 & $S_{11}$ \\ \hline
$S_{11}$ & $x^4+y^2z+xz^2$ & 16 & 4 5 6 & 2 5 6 & 3 4 4 & 11 & $-8$
 & $W_{13}$ \\ \hline
$S_{12}$ & $x^3y+y^2z+xz^2$ & 13 & 3 4 5 & 3 4 5 & 3 4 5 & 12 & 13
 & $S_{12}$ \\ \hline
$U_{12}$ & $x^4+y^3+z^3$ & 12 & 3 4 4 & 4 4 4 & 4 4 4 & 12 & 16
 & $U_{12}$ \\ \hline
\end{tabular}
\end{table}

\begin{figure}\centering
\unitlength1cm
\begin{picture}(8.5,7.5)
\put(0.5,0.5){\includegraphics{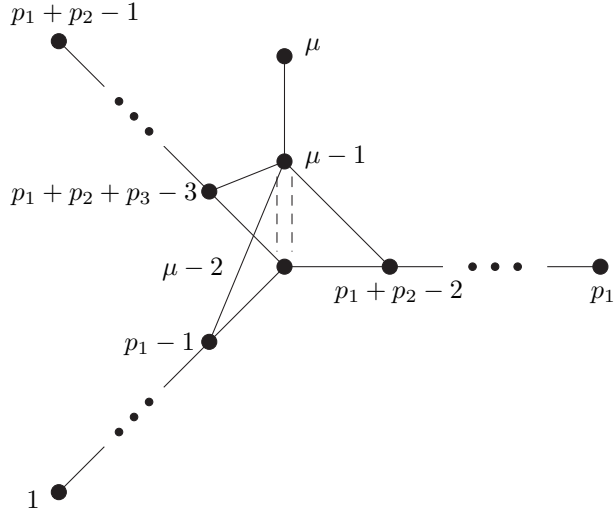}}
\put(0.2,0.4){$1$}
\put(1.5,2.5){$p_1-1$}
\put(2,3.5){$\mu-2$}
\put(0,4.5){$p_1+p_2+p_3-3$}
\put(0,6.9){$p_1+p_2-1$}
\put(4.3,3.2){$p_1+p_2-2$}
\put(7.7,3.2){$p_1$}
\put(3.9,5){$\mu-1$}
\put(3.9,6.5){$\mu$}
\end{picture}
\caption{Coxeter-Dynkin diagram of an exceptional unimodal singularity}
\label{Fig1}
\end{figure}

Let $(X,x_0)$ be one of the 14 hypersurface triangle singularities, and denote
by $X_t$ and $\mu$ its Milnor fibre and Milnor number respectively.
We denote by $\langle \ , \ \rangle$ the intersection form on $H_2(X_t,\ZZ)$
and by $H=(H_2(X_t,\ZZ),\langle \ , \ \rangle)$ the Milnor lattice.
A.~M.~Gabrielov
\cite{Gabrielov74} has shown that there exists a weakly distinguished basis of
vanishing cycles of
$H$ with a Coxeter-Dynkin diagram of the form of Fig.~\ref{Fig1}. The author
\cite{Ebeling81} has shown that this diagram even corresponds to a distinguished
basis of vanishing cycles (cf.\ also \cite{Ebeling96}). (For the notions of a
distinguished and weakly distinguished basis of vanishing cycles see e.g.\
\cite{AGV88}). The numbers $p_1$, $p_2$, $p_3$ are called the {\em Gabrielov
numbers} ${\rm Gab}(X)$ of the singularity. Here each vertex represents a
sphere of self-intersection number $-2$, two vertices connected by a
single solid edge have intersection number 1, and two vertices connected by a
double broken line have intersection number $-2$.
Using the results of K.~Saito (see \cite[Theorem~3.4.3]{Ebeling87}), one can
see that the Gabrielov numbers are uniquely determined by the singularity.
We denote by $d$ the discriminant of $H$, i.e.\ the determinant of an
intersection matrix with respect to a basis of $H$.

Arnold has now observed: There exists an involution $X \mapsto X^\ast$ on the
set of the 14 exceptional unimodal singularities, such that
$${\rm Dol}(X) = {\rm Gab}(X^\ast), \quad {\rm Gab}(X)={\rm Dol}(X^\ast), \quad
N=N^\ast, \quad \mu + \mu^\ast = 24. $$
This is called {\em Arnold's strange duality}. Note that also $d=d^\ast$.

H.~Pinkham \cite{Pinkham77} has given the following interpretation of this
duality. (This was independently also obtained by I.~V.~Dolgachev and
V.~V.~Nikulin \cite{DN77, Dolgachev82}.) The Milnor fibre $X_t$ can be
compactified in a weighted projective space to a surface with three cyclic
quotient singularities on the curve at infinity; a minimal resolution of these
singularities yields a K3 surface $S$. Denote by
$G(p_1,p_2,p_3)$ the subgraph of the graph of Fig.~\ref{Fig1} which is obtained
by omitting the vertices with indices
$\mu -1$ and $\mu$. Let $M(p_1,p_2,p_3)$ be the lattice (the free abelian group
with an integral quadratic form) determined by the graph $G(p_1,p_2,p_3)$. Then
$H = M(p_1,p_2,p_3)
\oplus U$, where $U$ is a unimodular hyperbolic plane (the lattice of rank 2
with a basis $\{e,e'\}$ such that $\langle e, e'\rangle = 1$, $\langle e,e
\rangle = \langle e',e' \rangle =0$) and
$\oplus$ denotes the orthogonal direct sum. The dual graph of the curve
configuration of
$S$ at infinity is given by $G(b_1,b_2,b_3)$. The inclusion $X_t \subset S$
induces a primitive embedding $H_2(X_t,\ZZ) \hookrightarrow H_2(S,\ZZ)$ and the
orthogonal complement is just the lattice $M(b_1,b_2,b_3)$. By \cite{Nikulin79}
the primitive embedding of $M(p_1,p_2,p_3) \oplus U$ into the unimodular K3
lattice $L:=H_2(S,\ZZ)$ is unique up to isomorphism.

In this way, Arnold's strange duality corresponds to a duality of K3 surfaces.
This is a two-dimensional analogue of the mirror symmetry between Calabi-Yau
threefolds. This has recently been worked out by Dolgachev \cite{Dolgachev95}.
We give an outline of his construction. Let $M$ be an even non-degenerate
lattice of signature $(1,t)$. An {\em $M$-polarized} K3 surface is a pair
$(S,j)$ where $S$ is a K3 surface and $j: M \hookrightarrow \mbox{Pic}(S)$ is a
primitive lattice embedding. Here $\mbox{Pic}(S)$ denotes the Picard group of
$S$. An $M$-polarized K3 surface $(S,j)$ is called {\em pseudo-ample} if $j(M)$
contains a pseudo-ample divisor class. We assume that $M$ has a unique
embedding into the K3 lattice $L$ and the orthogonal complement $M^\perp$ admits
an orthogonal splitting $M^\perp = U \oplus \check{M}$. (Dolgachev's
construction
is slightly more general.) Then we consider the complete family $\cal F$ of
pseudo-ample
$M$-polarized K3 surfaces and define its {\em mirror family} ${\cal F}^\ast$ to
be any complete family of pseudo-ample $\check{M}$-polarized K3 surfaces. It is
shown in \cite{Dolgachev95} that this is well defined and that there
is the following relation between $\cal F$ and ${\cal F}^\ast$: The dimension
of the family $\cal F$ is equal to the rank of the Picard group of a general
member from the mirror family ${\cal F}^\ast$. In particular, this can be
applied to $M=M(b_1,b_2,b_3)$ and $\check{M}= M(p_1,p_2,p_3)$ for one of the
14 Dolgachev triples $(b_1,b_2,b_3)$. See \cite{Dolgachev95} for further
results and references.

It was observed by I.~Nakamura \cite{Nakamura80, Nakamura81} that
Arnold's strange duality corresponds to Hirzebruch-Zagier duality of hyperbolic
(alias cusp) singularities. For details see \cite{Nakamura80, Nakamura81, EW85}.

\section{Kobayashi's duality of weight systems}

In his paper \cite{Kobayashi95}, M.~Kobayashi has observed a new feature of
Arnold's strange duality which we now want to explain.

A quadruple $W = (w_1, w_2, w_3; N)$ of positive integers with $N \in
\NN w_1+\NN w_2+\NN w_3$ is called a {\em weight system}. The integers $w_i$ are
called the weights and $N$ is called the degree of $W$. A weight system $W =
(w_1, w_2, w_3; N)$ is called {\em reduced} if $\gcd(w_1,w_2,w_3) = 1$.

Let $W = (w_1, w_2, w_3; N)$ and $W' = (w'_1, w'_2, w'_3; N')$ be two reduced
weight systems. An $3 \times 3$- matrix $Q$ whose elements are non-negative
integers is called a {\em weighted magic square} for $(W,W')$, if
$$(w_1,w_2,w_3)Q=(N,N,N) \quad \mbox{and} \quad
Q \left( \begin{array}{c} w'_1 \\ w'_2 \\ w'_3 \end{array} \right) =
\left( \begin{array}{c} N' \\ N' \\ N' \end{array} \right).$$
(In the case $w_1=w_2=w_3=w'_1=w'_2=w'_3=1$, $Q$ is an ordinary magic square.)
$Q$ is called {\em primitive}, if $| \det Q | = N = N'$. We say that the
weight systems $W$ and $W'$ are {\em dual} if there exists a primitive
weighted magic square for $(W,W')$.

Kobayashi now proves:

\begin{theorem}[M.~Kobayashi] \label{thm:Kobayashi}
Let $W = (w_1, w_2, w_3; N)$ be the weight system of one of the 14 exceptional
unimodal singularities. Then there exists a unique (up to permutation) dual
weight system $W^\ast$. The weight system $W^\ast$ belongs to the dual
singularity in the sense of Arnold.
\end{theorem}

Moreover, Kobayashi shows that there is a relation between this duality of
weight systems and the polar duality between certain polytopes associated to
the weight systems. Such a polar duality was considered by V.~Batyrev
\cite{Batyrev94} in connection with the mirror symmetry of Calabi-Yau
hypersurfaces in toric varieties. We refer to \cite{Ebeling98} for a more
precise discussion of this relation.

\section{Saito's duality of characteristic polynomials}

Let $f: (\CC^3,0) \to (\CC,0)$ be a germ of an analytic function defining an
isolated hypersurface singularity $(X,x_0)$. A characteristic homeomorphism of
the Milnor fibration of $f$ induces an automorphism $c:H_2(X_t,\ZZ) \to
H_2(X_t,\ZZ)$ called the {\em (classical) monodromy operator} of $(X,x_0)$. It
is a well known theorem (see e.g. \cite{Brieskorn70})  that the eigenvalues of
$c$ are roots of unity. This means that the characteristic polynomial
$\phi(\lambda) = \det (\lambda I - c)$ of $c$ is a monic polynomial the roots of
which are roots of unity. Such a polynomial can be written in the form
$$\phi(\lambda)= \prod_{m \geq 0} (\lambda^m -1)^{\chi_m} \quad \mbox{for} \
\chi_m \in \ZZ,$$
where all but finitely many of the integers $\chi_m$ are equal to zero. We note
some useful formulae.

\begin{proposition} \label{formulae}
\begin{itemize}
\item[{\rm (i)}] $\mu=\deg \phi = \sum_{m > 0} m\chi_m.$
\item[{\rm (ii)}] If $\sum_{m > 0} \chi_m = 0$ then
$$\phi(1) = \prod_{m > 0} m^{\chi_m}.$$
\item[{\rm (iii)}] $\phi(1) = (-1)^\mu d$.
\item[{\rm (iv)}] $ {\rm tr}\, c^k = \sum_{m | k} m \chi_m.$ \\
In particular $ {\rm tr}\, c = \chi_1$.
\end{itemize}
\end{proposition}

\noindent {\em Proof.} (i) is obvious. For the proof of (ii) we use the identity
$$(\lambda^m -1)=(\lambda -1)(\lambda^{m-1} + \ldots + \lambda +1).$$
To prove (iii), let $A$ be the intersection matrix with respect to a
distinguished basis $\{\delta_1, \ldots , \delta_\mu \}$ of vanishing cycles.
Write $A$ in the form $A=V+V^t$ where $V$ is an upper triangular matrix with
$-1$ on the diagonal. Let $C$ be the matrix of $c$ with respect to $\{\delta_1,
\ldots , \delta_\mu \}$. Then $C = - V^{-1}V^t$ (see e.g.\
\cite[Proposition~1.6.3]{Ebeling87}). Therefore
$$\phi(1)=\det (I-C)=\det V^{-1}(V+V^t) = (-1)^\mu d.$$
Finally, (iv) is obtained as in \cite{A'Campo75} using the identity
$$\det (I-tc) = \exp (\mbox{tr}\, (\log (I-tc))) = \exp (- \sum_{k \geq 1}
\frac{t^k}{k} \mbox{tr}\, c^k).$$
This proves Proposition~\ref{formulae}.

\addvspace{3mm}

By A'Campo's theorem
\cite{A'Campo73}
$$ \mbox{tr}\, c = -1.$$
We assume that $c$ has finite order $h$. This is e.g.\ true if $f$ is a
weighted homogeneous polynomial of degree $N$. In this case $h=N$. Then
$\chi_m=0$ for all
$m$ which do not divide $h$. K.~Saito
\cite{Saito94} defines a {\em dual polynomial}
$\phi^\ast(\lambda)$ to $\phi(\lambda)$:
$$\phi^\ast(\lambda) = \prod_{k | h} (\lambda^k -1)^{-\chi_{h/k}}.$$
He obtains the following result.

\begin{theorem}[K.~Saito] \label{thm:Saito}
If $\phi(\lambda)$ is the characteristic polynomial of the monodromy of an
 exceptional unimodal singularity $X$, then $\phi^\ast(\lambda)$ is the
corresponding polynomial of the dual singularity $X^\ast$.
\end{theorem}

For $\phi(\lambda)= \prod_{m|h} (\lambda^m -1)^{\chi_m}$ we use the symbolic
notation
$$\pi:= \prod_{m|h} m^{\chi_m}.$$
In the theory of finite groups, this symbol is known as a {\em Frame shape}
\cite{Frame, CN79}. For, if one has a rational finite-dimensional representation
of a finite group, then the zeros of the characteristic polynomials of each
element of the group are also roots of unity. For a given rational
representation, one can thus assign to each conjugacy class of the group its
Frame shape. The number
$$ \deg (\pi) = \sum m \chi_m$$
is called the {\em degree} of the Frame shape $\pi$.

Let us denote the Frame shape of the dual polynomial $\phi^\ast(\lambda)$ by
$\pi^\ast$. The Frame shapes of the monodromy operators of the 14 exceptional
unimodal singularities are listed in Table~\ref{Table2}.

\begin{table}\centering
\caption{Frame shapes of the 14 exceptional unimodal singularities}
\label{Table2}
\begin{tabular}{|c|c|c|c|}  \hline
Name   & $\pi$ & $\pi^\ast$ & Dual
\\ \hline
$E_{12}$ & $2 \cdot 3 \cdot 7 \cdot 42 / 1 \cdot 6 \cdot 14 \cdot 21$
 & $2 \cdot 3 \cdot 7 \cdot 42 / 1 \cdot 6 \cdot 14 \cdot 21$
 & $E_{12}$
\\ \hline
$E_{13}$ & $2 \cdot 3 \cdot 30 / 1 \cdot 6 \cdot 15$ & $2 \cdot 5 \cdot 30
/ 1 \cdot 10 \cdot 15$ & $Z_{11}$  \\ \hline
$E_{14}$ & $2 \cdot 3 \cdot 24 / 1 \cdot 6 \cdot 8$  & $3 \cdot 4 \cdot 24
/ 1 \cdot 8 \cdot 12$ & $Q_{10}$  \\ \hline
$Z_{12}$ & $2 \cdot 22 / 1 \cdot 11$ & $2 \cdot 22 / 1 \cdot
11$ & $Z_{12}$  \\ \hline
$Z_{13}$ & $2 \cdot 18 / 1 \cdot 6$  & $3 \cdot 18 / 1 \cdot 9$ & $Q_{11}$
 \\ \hline
$Q_{12}$ & $3 \cdot 15 / 1 \cdot 5$ & $3 \cdot 15 / 1 \cdot 5$ & $Q_{12}$
 \\ \hline
$W_{12}$ & $2 \cdot 5 \cdot 20 / 1 \cdot 4 \cdot 10$ & $2 \cdot 5 \cdot 20
/ 1 \cdot 4 \cdot 10$ & $W_{12}$  \\
\hline
$W_{13}$ & $2 \cdot 16 / 1 \cdot 4$ & $4 \cdot 16 / 1 \cdot 8$ & $S_{11}$
 \\ \hline
$S_{12}$ & $13 / 1$ & $13 / 1$ & $S_{12}$  \\
\hline
$U_{12}$ & $4 \cdot 12 / 1 \cdot 3$ & $4 \cdot 12 / 1 \cdot 3$ & $U_{12}$
 \\ \hline
\end{tabular}
\end{table}

To two Frame shapes $\pi = \prod m^{\chi_m}$ and $\pi'=\prod m^{\chi'_m}$ of
degree $n$ and $n'$ respectively one can associate a Frame shape $\pi\pi'$ of
degree $n+n'$ by concatenation
$$\pi\pi':=\prod m^{\chi_m}\prod m^{\chi'_m}=\prod m^{\chi_m + \chi'_m}.$$
If $\pi$ and $\pi'$ are the Frame shapes of the operators $c: H \to H$ and $c'
: H' \to H'$ respectively, $\pi\pi'$ is the Frame shape of the operator $c
\oplus c' : H \oplus H' \to H \oplus H'$.

In the appendix of
\cite{Saito94}, Saito notes the following observation: If $\pi$ is the Frame
shape of the monodromy operator of an exceptional unimodal singularity, then the
symbol
$\pi\pi^\ast$ of degree 24 is a Frame shape of a conjugacy class of the
automorphism group of the Leech lattice. The Leech lattice is a 24-dimensional
even unimodular positive definite lattice which contains no roots (see e.g.
\cite{Ebeling94}). It was discovered by J.~Leech in connection with the search
for densest sphere packings. Its automorphism group $G$ is a group of order
$2^{22}3^95^47^211\cdot 13 \cdot 23$. The quotient group
$\mbox{Co}_1 := G/\{\pm 1\}$ is a famous sporadic simple group discovered by
J.~Conway. The Frame shapes of the 164 conjugacy classes of $G$ have been
listed by T.~Kondo \cite{Kondo85}.

\section{An extension of Arnold's strange duality}

C.~T.~C.~Wall and the author \cite{EW85} have found an extension of Arnold's
strange duality. In order to consider this extension, we have to enlarge the
class of singularities which we want to discuss.

On the one hand, instead of restricting to the hypersurface case, we can also
look at isolated complete intersection singularities (abbreviated ICIS in the
sequel). Pinkham has shown \cite{Pinkham77b} that for exactly 8 triples
$(b_1,b_2,b_3)$ the triangle singularities with these Dolgachev numbers are
ICIS, but not hypersurface singularities. They are given by germs of
analytic mappings
$(g,f): (\CC^4,0) \to (\CC^2, 0)$. They are
$\cal K$-unimodal singularities and appear in Wall's classification
\cite{Wall83}. For certain values of the module, the equations are again
weighted homogeneous. The corresponding 8 triples $(b_1,b_2,b_3)$, Wall's names
and weighted homogeneous equations with weights $(w_1,w_2,w_3,w_4)$ and degrees
$(N_1,N_2)$ are indicated in Table~\ref{Table3}.

\begin{table}\centering
\caption{The 8 triangle ICIS}
\label{Table3}
{\small
\begin{tabular}{|c|c|c|c|c|c|c|c|c|}  \hline
Name   & Equations & $N$ & Weights & Dol & Gab & $\mu$ & $d$  &
Dual
\\ \hline
$J'_9$ & $\begin{array}{c} xw+y^2 \\ x^3+yw+z^2 \end{array}$
 &       $\begin{array}{c} 16 \\ 18 \end{array}$
 & 6 8 9 10 & 2 3 10 & 2 2 2 3 & 9 & $-4$
 & $J_{3,0}$ \\ \hline
$J'_{10}$ & $\begin{array}{c} xw+y^2 \\ x^2y+yw+z^2 \end{array}$
 &          $\begin{array}{c} 12 \\ 14 \end{array}$
 & 4 6 7 8 & 2 4 8 & 2 2 2 4 & 10 & 8
 & $Z_{1,0}$ \\ \hline
$J'_{11}$ & $\begin{array}{c} xw+y^2 \\ x^4+yw+z^2 \end{array}$
 &          $\begin{array}{c} 10 \\ 12 \end{array}$
 & 3 5 6 7 & 3 3 7 & 2 2 2 5 & 11 & $-12$
 & $Q_{2,0}$ \\ \hline
$K'_{10}$ & $\begin{array}{c} xw+y^2 \\ x^3+z^2+w^2 \end{array}$
 &          $\begin{array}{c} 10 \\ 12 \end{array}$
 & 4 5 6 6 & 2 6 6 & 2 3 2 3 & 10 & 12
 & $W_{1,0}$ \\ \hline
$K'_{11}$ & $\begin{array}{c} xw+y^2 \\ x^2y+z^2+w^2 \end{array}$
 &          $\begin{array}{c} 8 \\ 10 \end{array}$
 & 3 4 5 5 & 3 5 5 & 2 3 2 4 & 11 & $-20$
 & $S_{1,0}$ \\ \hline
$L_{10}$ & $\begin{array}{c} xw+yz \\ x^3+yw+z^2 \end{array}$
 &         $\begin{array}{c} 11 \\ 12 \end{array}$
 & 4 5 6 7 & 2 5 7 & 2 2 3 3 & 10 & 11
 & $W_{1,0}$ \\ \hline
$L_{11}$ & $\begin{array}{c} xw+yz \\ x^2y+yw+z^2 \end{array}$
 &         $\begin{array}{c} 9 \\ 10 \end{array}$
 & 3 4 5 6 & 3 4 6 & 2 2 3 4 & 11 & $-18$
 & $S_{1,0}$ \\ \hline
$M_{11}$ & $\begin{array}{c} 2xw+y^2+z^2 \\ x^3+2yw \end{array}$
 &         $\begin{array}{c} 8 \\ 9 \end{array}$
 & 3 4 4 5 & 4 4 5 & 2 3 3 3 & 11 & $-24$
 & $U_{1,0}$ \\ \hline
\end{tabular}}
\end{table}

\begin{figure}\centering
\unitlength1cm
\begin{picture}(8.5,8)
\put(0.5,0.5){\includegraphics{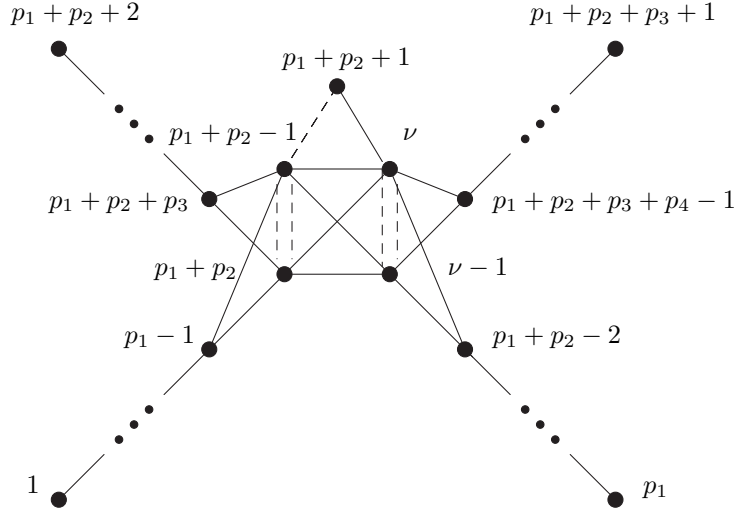}}
\put(0.2,0.7){$1$}
\put(8.4,0.7){$p_1$}
\put(1.5,2.7){$p_1-1$}
\put(6.4,2.7){$p_1+p_2-2$}
\put(1.9,3.6){$p_1+p_2$}
\put(5.8,3.6){$\nu-1$}
\put(0.5,4.5){$p_1+p_2+p_3$}
\put(2.1,5.4){$p_1+p_2-1$}
\put(5.2,5.4){$\nu$}
\put(6.4,4.5){$p_1+p_2+p_3+p_4-1$}
\put(0,7){$p_1+p_2+2$}
\put(3.6,6.4){$p_1+p_2+1$}
\put(6.9,7){$p_1+p_2+p_3+1$}
\end{picture}
\caption{Coxeter-Dynkin diagram of a triangle ICIS}
\label{Fig2}
\end{figure}

By \cite{Hamm71}, the notion of Milnor fibre can also be extended to ICIS. We
assume that $(g,f)$ are generically chosen such that $(X',0)=(g^{-1}(0),0)$
is an
isolated hypersurface singularity of minimal Milnor number $\mu_1$ among such
choices of $g$. Then the {\em monodromy operator} of $(X,0)$ is defined to be
the monodromy operator of the function germ $f: (X',0) \to (\CC,0)$. By
\cite{Ebeling87} there exists a distinguished set of generators consisting of
$\nu:=\mu + \mu_1$ vanishing cycles, where $\mu$ is the rank of the second
homology group of the Milnor fibre. Again the monodromy operator is the Coxeter
element of this set, i.e.\ the product of the
$\nu$ reflections corresponding to the vanishing cycles of the distinguished set
of generators. For the 8 triangle ICIS, a Coxeter-Dynkin diagram
corresponding to
such a distinguished set is depicted in Fig.~\ref{Fig2} (cf.\ \cite{Ebeling87}).
Let us call the characteristic numbers $p_1$, $p_2$, $p_3$, $p_4$ of these
graphs
the {\em Gabrielov numbers} ${\rm Gab}(X)$ of the singularity. They are also
indicated in Table~\ref{Table3}. Again, using \cite[Theorem~3.4.3]{Ebeling87}
one can see that these numbers are uniquely defined.

On the other hand, instead of starting with a hyperbolic triangle, we can start
with a hyperbolic quadrilateral. Let $b_1$, $b_2$, $b_3$, $b_4$ be positive
integers such that
$$\frac{1}{b_1} + \frac{1}{b_2} + \frac{1}{b_3} + \frac{1}{b_4} < 2.$$
One can perform the same construction as above with a solid quadrilateral with
angles $\frac{\pi}{b_1}$, $\frac{\pi}{b_2}$,
$\frac{\pi}{b_3}$, $\frac{\pi}{b_4}$ instead of a triangle. The
resulting normal surface singularities are called {\em quadrilateral
singularities} (with {\em Dolgachev numbers} $b_1$, $b_2$, $b_3$, $b_4$) (cf.\
\cite{Looijenga84}). Again these singularities admit a natural good
$\CC^\ast$-action. A minimal good resolution consists of a rational curve of
self-intersection number
$-2$ together with four rational curves of self-intersection numbers $-b_1$,
$-b_2$,
$-b_3$, and $-b_4$ respectively intersecting the first curve transversely.

For 6 quadruples
$(b_1,b_2,b_3,b_4)$ these singularities are isolated hypersurface singularities.
They are bimodal in the sense of Arnold. They can again be defined by weighted
homogeneous equations. By
\cite{EW85} they have a distinguished basis with a Coxeter-Dynkin diagram of the
following form: It is obtained by adding a new vertex with number 0 to the graph
of Fig.~\ref{Fig1} in one of the following two ways:
\begin{itemize}
\item[(1)] It is connected to the vertex $\mu$  and to the vertex
$p_1+p_2$ by ordinary lines. We refer to this diagram by the symbol $(p_1,
p_2, \underline{p_3})$.
\item[(2)] It is connected to the vertex $\mu$ and to the vertices
$p_1$ and $p_1+p_2-1$ by ordinary lines. We refer to this diagram by the symbol
$(p_1, \underline{p_2}, \underline{p_3})$.
\end{itemize}
The corresponding data are indicated in Table~\ref{Table4}. The numbers $p_1$,
$p_2$, $p_3$ defined by this procedure are unique except for the singularities
$W_{1,0}$ and $S_{1,0}$, where we have two triples each. This follows from
\cite[Supplement to Theorem~4.1]{Ebeling83}. One can also obtain an easier
proof by considering the discriminants (and in one case the discriminant
quadratic form) of the singularities and the possible determinants of the
graphs (and in one case the discriminant quadratic form determined by the
graph). The absolute values of the determinants in the respective cases are
\begin{itemize}
\item[(1)] $|d|=4(p_1 p_2 -p_1 - p_2)$,
\item[(2)] $|d|=(p_1 -1)(p_2+p_3)$.
\end{itemize}

\begin{table}\centering
\caption{The 6 quadrilateral hypersurface singularities}
\label{Table4}
{\small
\begin{tabular}{|c|c|c|c|c|c|c|c|c|}  \hline
Name   & Equation & $N$ & Weights & Dol & Gab & $\mu$ & $d$ &
Dual
\\ \hline
$J_{3,0}$ & $x^9+y^3+z^2$ & 18 & 2 6 9 & 2 2 2 3 & 2 3 \underline{10}
 & 16 & 4
 & $J'_9$ \\ \hline
$Z_{1,0}$ & $x^7+xy^3+z^2$ & 14 & 2 4 7 & 2 2 2 4 & 2 4 \underline{8}
 & 15 & $-8$
 & $J'_{10}$ \\ \hline
$Q_{2,0}$ & $x^6+y^3+xz^2$ & 12 & 2 4 5 & 2 2 2 5 & 3 3 \underline{7}
 & 14 & 12
 & $J'_{11}$ \\ \hline
$W_{1,0}$ & $x^6+y^4+z^2$ & 12 & 2 3 6 & 2 2 3 3
 & $\begin{array}{c}
    2 \  5 \ \underline{7} \\ 2 \ \underline{6} \ \underline{6}
    \end{array}$
 & 15 & $-12$
 & $\begin{array}{c} K'_{10} \\ L_{10} \end{array}$ \\ \hline
$S_{1,0}$ & $x^5+y^2z+xz^2$ & 10 & 2 3 4 & 2 2 3 4
 & $\begin{array}{c}
    3 \ 4 \ \underline{6} \\ 3 \ \underline{5} \ \underline{5}
    \end{array}$
 & 14 & 20
 & $\begin{array}{c} K'_{11} \\ L_{11} \end{array}$ \\ \hline
$U_{1,0}$ & $x^3y+y^3+z^3$ & 9 & 2 3 3 & 2 3 3 3 & 4 \underline{4}
 \underline{5} & 14 & 27
 & $M_{11}$ \\ \hline
\end{tabular}}
\end{table}

For another 5 quadruples $(b_1,b_2,b_3,b_4)$ the quadrilateral singularities
with these Dolgachev numbers are ICIS. They are also $\cal K$-unimodal and
appear in Wall's lists \cite{Wall83}. They can also be given by weighted
homogeneous
equations. These equations and Wall's names are listed in Table~\ref{Table5}.

\begin{table}\centering
\caption{The 5 quadrilateral ICIS}
\label{Table5}
\begin{tabular}{|c|c|c|c|c|}  \hline
Name   & Equations & Restrictions & $N$ & Weights
\\ \hline
$J'_{2,0}$ & $\begin{array}{c} xw+y^2 \\ ax^5+xy^2+yw+z^2
              \end{array}$ & $a \neq 0,-\frac{4}{27}$
 &           $\begin{array}{c} 8 \\ 10 \end{array}$
 & 2 4 5 6  \\ \hline
$L_{1,0}$ & $\begin{array}{c} xw+yz \\ ax^4+xy^2+yw+z^2
              \end{array}$ & $a \neq 0,-1$
 &          $\begin{array}{c} 7 \\ 8 \end{array}$
 & 2 3 4 5  \\ \hline
$K'_{1,0}$ & $\begin{array}{c} xw+y^2 \\ ax^4+xy^2+z^2+w^2
              \end{array}$ & $a \neq 0,\frac{1}{4}$
 &           $\begin{array}{c} 6 \\ 8 \end{array}$
 & 2 3 4 4  \\ \hline
$M_{1,0}$ & $\begin{array}{c} 2xw+y^2-z^2 \\ x^2y+ax^2z+2yw
              \end{array}$ & $a \neq \pm 1$
 &          $\begin{array}{c} 6 \\ 7 \end{array}$
 & 2 3 3 4 \\ \hline
$I_{1,0}$ & $\begin{array}{c} x^3+w(y-z) \\ ax^3+y(z-w)
              \end{array}$ & $a \neq 0,1$
 &          $\begin{array}{c} 6 \\ 6 \end{array}$
 & 2 3 3 3  \\ \hline
\end{tabular}
\begin{tabular}{|c|c|c|c|c|c|}  \hline
Name    & Dol & Gab & $\mu$ & $d$  & Dual
\\ \hline
$J'_{2,0}$  & 2 2 2 6 & 2 2 2 \underline{6} & 13 & $-16$
 & $J'_{2,0}$ \\ \hline
$L_{1,0}$ &  2 2 3 5 & $\begin{array}{c} \mbox{2 2 3 \underline{5}} \\
\mbox{2 2 \underline{4} \underline{4}} \end{array}$ & 13 & $-28$
 & $\begin{array}{c} L_{1,0} \\ K'_{1,0} \end{array}$ \\ \hline
$K'_{1,0}$ &  2 2 4 4 & $\begin{array}{c} \mbox{2 3 2 \underline{5}} \\
\mbox{2 \underline{4} 2 \underline{4}}  \end{array}$  & 13 & $-32$
 & $\begin{array}{c} L_{1,0} \\ K'_{1,0} \end{array}$ \\ \hline
$M_{1,0}$ &  2 3 3 4 & $\begin{array}{c} \mbox{2 3 \underline{3} \underline{4}}
\\ \mbox{2 \underline{3} 3 \underline{4}} \end{array}$ & 13 &
$-42$
 & $M_{1,0}$ \\ \hline
$I_{1,0}$ &  3 3 3 3 &
$\underline{3}$ $\underline{3}$ $\underline{3}$ $\underline{3}$
 & 13
 & $-54$
 & $I_{1,0}$ \\ \hline
\end{tabular}
\end{table}

The singularities $J'_{2,0}$, $L_{1,0}$, $K'_{1,0}$, and $M_{1,0}$ can be given
by equations $(g,f)$ where $g$ has Milnor number $\mu_1 = 1$. Coxeter-Dynkin
diagrams of these singularities are computed in \cite{Ebeling87}. Using
transformations as in the proof of \cite[Proposition~3.6.1]{Ebeling87}, these
graphs can be transformed to the following graphs. A Coxeter-Dynkin diagram
corresponding to  a distinguished set of generators is obtained by adding a
new vertex to the graph of Fig.~\ref{Fig2}. It gets the number $p_1+p_2+2$ and
the indices of the old vertices with numbers $p_1+p_2+2, p_1+p_2+3, \ldots ,
\nu$ are shifted by $1$. New edges are introduced in one of the following ways:
\begin{itemize}
\item[(1)] The new vertex is connected to the vertex $p_1+p_2+1$ and
to the vertex with new index $p_1+p_2+p_3+3$ (old index $p_1+p_2+p_3+2$) by
ordinary lines. We refer to this diagram by the symbol $(p_1, p_2, p_3,
\underline{p_4})$.
\item[(2)] The new vertex is connected to the vertex $p_1+p_2+1$
and to the vertices with new indices $p_1+p_2+3$ and $p_1+p_2+p_3+2$ (old
indices $p_1+p_2+2$ and $p_1+p_2+p_3+1$ respectively) by ordinary lines. We
refer to this diagram by the symbol $(p_1,p_2,\underline{p_3},\underline{p_4})$.
\item[(3)] The new vertex is connected to the vertex $p_1+p_2+1$, to the
vertex $p_1$, and to the vertex with new index $p_1+p_2+p_3+2$ (old index
$p_1+p_2+p_3+1$) by ordinary lines. We refer to this diagram by the symbol
$(p_1,\underline{p_2},p_3,\underline{p_4})$.
\end{itemize}
The absolute values of the determinants of the respective graphs are
\begin{itemize}
\item[(1)] $|d|=4(p_1p_2p_3 -p_1-p_3)$,
\item[(2)] $|d|=p_1p_2(p_3+p_4+2)-p_1-p_2-p_3-p_4$,
\item[(3)] $|d|=p_1p_3(p_2+p_4)$.
\end{itemize}
Comparing the values of these determinants with the discriminants of the above 4
quadrilateral ICIS, we find that the graphs listed in Table~\ref{Table5} are
the only possible graphs of the types (1), (2), or (3) for these singularities.
Again, in two cases the numbers $p_1$, $p_2$, $p_3$, $p_4$ are not
uniquely defined.

\begin{figure}\centering
\unitlength1cm
\begin{picture}(5.5,5.5)
\put(0.5,0.5){\includegraphics{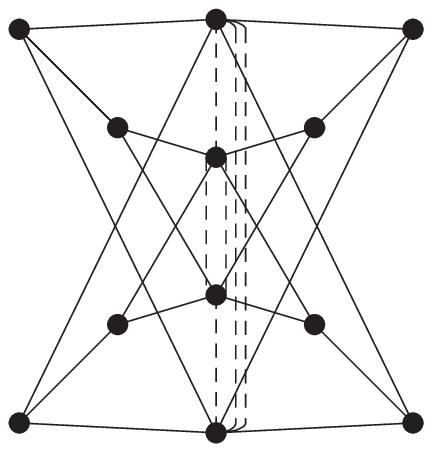}}
\put(2.3,3.6){$1$}
\put(1.6,3.9){$2$}
\put(1.6,1.3){$3$}
\put(3.5,3.9){$4$}
\put(3.5,1.3){$5$}
\put(2.3,1.6){$6$}
\put(2.5,0.1){$7$}
\put(0.2,4.6){$8$}
\put(0.2,0.5){$9$}
\put(4.9,4.6){$10$}
\put(4.9,0.5){$11$}
\put(2.5,5.1){$12$}
\end{picture}
\caption{Subgraph of a Coxeter-Dynkin diagram of the singularity $I_{1,0}$}
\label{Fig3}
\end{figure}

For the remaining singularity $I_{1,0}$, $\mu_1=2$. This singularity can
be given by the following equations:
\begin{eqnarray*}
g(z) & = & z_1^2 + z_2^2 + z_3^2 + z_4^3, \\
f(z) & = & a_1z_1^2 + a_2z_2^2 + a_3z_3^2 + a_4z_4^3,
\end{eqnarray*}
where $a_i \in \RR$, $a_1 < a_2 < a_3 < a_4$. For such a singularity H.~Hamm
\cite{Hamm72} has given a basis of the complexified Milnor lattice
$H_\CC = H
\otimes \CC$. As in \cite[Sect.~2.3]{Ebeling87}, one can show that the cycles
he constructs are in fact vanishing cycles and that there exists a
distinguished set $\{\delta_1, \ldots , \delta_\nu\}$ of generators for
this singularity with the following intersection numbers:
$$\langle \delta_{i+1},\delta_{i+2} \rangle = -1, \quad i=0,2,4,6,8,10,$$
$$\langle \delta_{i+1},\delta_{i+3} \rangle =\langle \delta_{i+1},\delta_{i+4}
\rangle = 0, \langle \delta_{i+2},\delta_{i+3} \rangle =\langle
\delta_{i+2},\delta_{i+4} \rangle = 0,\quad i=0,4,8,$$
$$\langle \delta_{i+1},\delta_{j} \rangle =\langle \delta_{i+3},\delta_{j}
\rangle = -1, \quad i=0,4,8, \ 1 \leq j\leq 12, j\neq i+1,i+2,i+3,i+4,$$
$$\langle \delta_{i+2},\delta_{2k} \rangle =\langle \delta_{i+4},\delta_{2k}
\rangle = -1, \quad i=0,4,8, \ 1 \leq k \leq 6, 2k \neq i+2, i+4,$$
$$\langle \delta_{i+2},\delta_{2k-1} \rangle =\langle \delta_{i+4},\delta_{2k-1}
\rangle = 0, \quad i=0,4,8, \ 1 \leq k \leq 6, 2k-1 \neq i+1, i+3,$$
$$\langle \delta_{2k-1} , \delta_{13} \rangle = -1, \langle \delta_{2k-1} ,
\delta_{15} \rangle = 0, \langle \delta_{2k} , \delta_{13} \rangle = 0, \langle
\delta_{2k} , \delta_{15} \rangle = -1,
 \quad 1 \leq k \leq 6,$$
$$\langle \delta_{i} , \delta_{14} \rangle = -1, \quad 1 \leq i \leq 12,$$
$$\langle \delta_{13} , \delta_{14} \rangle = 0, \langle \delta_{13} ,
\delta_{15} \rangle = 0, \langle \delta_{14} , \delta_{15} \rangle = 0.$$
By the following sequence of braid group transformations (for the notation see
e.g. \cite{Ebeling87}) the distinguished set $\{\delta_1,
\ldots , \delta_\nu\}$ can be transformed to a distinguished set $\{\delta'_1,
\ldots , \delta'_\nu\}$ where the Coxeter-Dynkin diagram corresponding to
the subset $\{\delta'_1, \ldots , \delta'_{12}\}$ is the graph of
Fig.~\ref{Fig3}:
$$\beta_5,\beta_4,\beta_3,\beta_2;\alpha_6,\alpha_7,\alpha_8,\alpha_9,\alpha
_{10},
\alpha_{11};\alpha_5,\alpha_6;\alpha_3,\alpha_4,\alpha_5;$$
$$\beta_8,\beta_7,\beta_6;\beta_{10},\beta_9,\beta_8,\beta_7;\kappa_4,\kappa_5,
\kappa_6,\kappa_7,\kappa_8,\kappa_9.$$
We refer to the Coxeter-Dynkin diagram corresponding to $\{\delta'_1,
\ldots , \delta'_\nu\}$
by the symbol $(\underline{3},\underline{3},\underline{3},\underline{3})$. This
notation is motivated by Fig.~3. We admit, however, that it is somewhat
arbitrary.

The corresponding symbols are indicated in Table~\ref{Table5}.

If one compares the Dolgachev and Gabrielov numbers of Tables \ref{Table3} and
\ref{Table4} and of Table~\ref{Table5}, then one observes a correspondence
between the 8 triangle ICIS and the 6 quadrilateral hypersurface singularities
and between the 5 quadrilateral ICIS. The corresponding ''dual'' singularities
are indicated in the last column of each table.
Note that this correspondence is not always a duality in the strict sense. For
the quadrilateral hypersurface singularity $W_{1,0}$ we have two corresponding
ICIS $K'_{10}$ and $L_{10}$ and for $S_{1,0}$ we have the corresponding ICIS
$K'_{11}$ and $L_{11}$. The quadrilateral ICIS $L_{1,0}$ and $K'_{1,0}$ are both
self-dual and dual to each other. In the other cases the correspondence
is one-to-one. Pinkham also defined "Gabrielov numbers" (in a weaker sense) for
the triangle ICIS
\cite{Pinkham77b} and he already made part of this observation (unpublished).
This duality also corresponds to the Hirzebruch-Zagier duality of
cusp singularities (see
\cite{Nakamura81, EW85}).

If one now compares the Milnor numbers of dual singularities,  one
finds
\begin{itemize}
\item for the triangle ICIS versus quadrilateral hypersurface singularities:
$\mu + \mu^\ast = 25$.
\item for the quadrilateral ICIS: $\mu + \mu^\ast = 26$.
\end{itemize}
(Note that also $d$ and $d^\ast$ do not coincide in each case.) So one still
has to alter something. There are two alternatives:
\begin{itemize}
\item[(1)] subtract 1 for quadrilateral.
\item[(2)] subtract 1 for ICIS.
\end{itemize}

In \cite{EW85} we considered the first alternative. The quadrilateral
singularities are first elements ($l=0$) of series of singularities indexed by
a non-negative integer $l$. We showed that to each such series one can
associate a virtual element $l=-1$.  We defined for these Milnor lattices,
Coxeter-Dynkin diagrams, and monodromy operators, and showed that all features
of Arnold's strange duality including Pinkham's interpretation continue to hold.
For more details see below.

\begin{figure}\centering
\unitlength1cm
\begin{picture}(7.5,7.5)
\put(0.5,0.5){\includegraphics{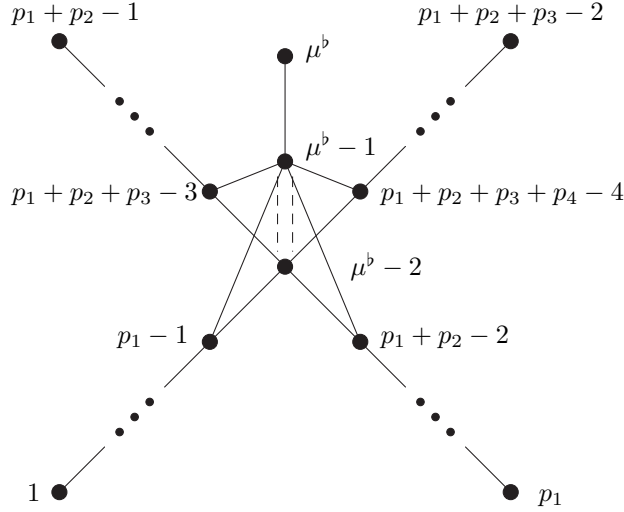}}
\put(0.2,0.5){$1$}
\put(7,0.5){$p_1$}
\put(1.4,2.6){$p_1-1$}
\put(4.9,2.6){$p_1+p_2-2$}
\put(4.5,3.5){$\mu^\flat-2$}
\put(0,4.5){$p_1+p_2+p_3-3$}
\put(4.9,4.5){$p_1+p_2+p_3+p_4-4$}
\put(0,6.9){$p_1+p_2-1$}
\put(5.4,6.9){$p_1+p_2+p_3-2$}
\put(3.9,5.1){$\mu^\flat-1$}
\put(3.9,6.4){$\mu^\flat$}
\end{picture}
\caption{Reduced Coxeter-Dynkin diagram $D^\flat$}
\label{Fig4}
\end{figure}

A new discovery is that the second alternative works as well, and this also
leads to an extension of Saito's duality. Recall that the triangle or
quadrilateral ICIS with $\mu_1 = 1$ have a Coxeter-Dynkin diagram $D$ which is
either the graph of Fig.~\ref{Fig2} or an extension of it. By similar
transformations as in the proof of \cite[Proposition~3.6.2]{Ebeling87}, this
graph can be transformed to a graph  containing the subgraph $D^\flat$ depicted
in Fig.~\ref{Fig4}.
(Unfortunately, this proof has to be modified slightly. The correct sequence of
transformations is
$$\beta_{\rho-1}, \beta_{\rho-2}, \beta_{\rho-3}, \beta_{\rho-3},
\beta_{\rho-2},
\kappa_{\rho-1}$$
and one has to consider the vertices
$\lambda^\prime_{\rho-2}$ and $\lambda^\prime_{\rho}$ instead of
$\lambda^\prime_{\rho-1}$ and $\lambda^\prime_{\rho}$.) By
\cite[Remark~3.6.5]{Ebeling87}, the passage from $D$ to $D^\flat$ can be
considered as a kind of "desuspension". For the singularity $I_{1,0}$, let
$D^\flat$ be the graph of Fig.~\ref{Fig3}. In each case, the new Coxeter-Dynkin
diagram
$D^\flat$ has
$\mu -1$ instead of $\mu + \mu_1$ vertices. Denote the corresponding Coxeter
element (product of reflections corresponding to the vanishing cycles) by
$c^\flat$. For the ICIS with $\mu_1=1$ one can compute (see
\cite[Proposition~3.6.2]{Ebeling87}) that
$$ \pi(c^\flat) = \pi(c)/1.$$
This means that $c^\flat$ has the same eigenvalues as $c$ but the multiplicities
of the eigenvalue 1 differ by 1. For the singularity $I_{1,0}$ one has
$\pi(c) = 6^3 / 1 \cdot 2^2$ (cf.\ \cite{Hamm72}), whereas $\pi(c^\flat) =
3^26^2 / 1^22^2$. Note that in any case
$$\mbox{tr}\, c^\flat = -2.$$
The passage from $c$ to
$c^\flat$ corresponds to the passage from the Milnor lattice $H$ of
rank
$\mu$
to a sublattice $H^\flat$ of rank $\mu^\flat = \mu -1$. The
corresponding discriminants and discriminant quadratic forms of the lattices
$H^\flat$ are listed in Table~\ref{Table6}. Here we use the notation of
\cite{EW85}.

\begin{table}\centering
\caption{Discriminants and discriminant quadratic forms of the lattices
$H^\flat$}
\label{Table6}
\begin{tabular}{|c|c|c|c|c|c|}  \hline
Name & $\mu^\flat$ & $d^\flat$ & $(H^\flat)^\ast/H^\flat$ & dual form & dual
\\ \hline
$J'_9$ & 8 & 4 & $q_{D_4}$ & $q_{D_4}$ & $J_{3,0}$ \\ \hline
$J'_{10}$ & 9 &$-8$ & $q_{D_4} + q_{A_1}$ & $q_{D_4} + w^1_{2,1}$ & $Z_{1,0}$
 \\ \hline
$J'_{11}$ & 10 & 12 & $q_{D_4} + q_{A_2}$ & $q_{D_4} + w^{-1}_{3,1}$ & $Q_{2,0}$
\\ \hline
$K'_{10}$ & 9 & $-12$ & $w^1_{2,2} + w^{-1}_{3,1}$
& $w^{-1}_{2,2} + w^1_{3,1}$ & $W_{1,0}$ \\ \hline
$L_{10}$ & 9 & $-12$ & $w^1_{2,2} + w^{-1}_{3,1}$
& $w^{-1}_{2,2} + w^1_{3,1}$ & $W_{1,0}$ \\ \hline
$K'_{11}$ & 10 & $20$ & $w^1_{2,1} + w^1_{2,1} + w^{-1}_{5,1}$
& $w^{-1}_{2,1} + w^{-1}_{2,1} + w^{-1}_{5,1}$ & $S_{1,0}$ \\ \hline
$L_{11}$ & 10 & $20$ & $w^1_{2,1} + w^1_{2,1} + w^{-1}_{5,1}$
& $w^{-1}_{2,1} + w^{-1}_{2,1} + w^{-1}_{5,1}$ & $S_{1,0}$ \\ \hline
$M_{11}$ & 10 & $27$ & $w^{-1}_{3,2} + w^{-1}_{3,1}$
& $w^1_{3,2} + w^1_{3,1}$ & $U_{1,0}$ \\ \hline
$J'_{2,0}$ & 12 & 16 & $v_1 + v_1$ & $v_1 + v_1$ & $J'_{2,0}$\\ \hline
$\begin{array}{c} L_{1,0} \\ K'_{1,0} \end{array}$ & 12 & 32 & $w^{1}_{2,2} +
w^{-1}_{2,3}$ &
$w^{1}_{2,2} + w^{-1}_{2,3}$  & $\begin{array}{c} L_{1,0} \\ K'_{1,0}
\end{array}$ \\
\hline
$M_{1,0}$ & 12 & 49 & $w^1_{7,1} + w^{-1}_{7,1}$
& $w^1_{7,1} + w^{-1}_{7,1}$ & $M_{1,0}$ \\ \hline
$I_{1,0}$ & 12 & 81 & $\begin{array}{l} w^1_{3,1} + w^1_{3,1} \\
 + w^{-1}_{3,1} + w^{-1}_{3,1} \end{array}$   &
$\begin{array}{l} w^{-1}_{3,1} + w^{-1}_{3,1} \\
+ w^1_{3,1} + w^1_{3,1} \end{array}$   & $I_{1,0}$ \\
\hline
\end{tabular}
\end{table}

The Frame shapes of the corresponding operators $c^\flat$ are listed in
Table~\ref{Table7}.

\begin{table}\centering
\caption{Frame shapes of the triangle ICIS and quadrilateral singularities}
\label{Table7}
\begin{tabular}{|c|c|c|c|}  \hline
Name   & $\pi$ & $\pi^\ast$ & Dual
\\ \hline
$J_{3,0}$ & $2 \cdot 3 \cdot 18^2 / 1 \cdot 6 \cdot 9^2$  & $2^23 \cdot 18
/ 1^26 \cdot 9$ & $J'_9$ \\ \hline
$Z_{1,0}$ & $2 \cdot 14^2 / 1 \cdot 7^2$  & $2^2 14 / 1^27$ &
$J'_{10}$  \\ \hline
$Q_{2,0}$ & $3 \cdot 12^2 / 1 \cdot 6^2$  & $2^2 12 / 1^24$ &
$J'_{11}$ \\ \hline
$W_{1,0}$ & $2 \cdot 12^2 / 1 \cdot 4 \cdot 6$  & $2 \cdot 3 \cdot 12 /
1^26$
 & $\begin{array}{c} K'_{10} \\ L_{10} \end{array}$  \\ \hline
$S_{1,0}$ & $10^2 / 1 \cdot 5$  & $2 \cdot 10 / 1^2$
 & $\begin{array}{c} K'_{11} \\ L_{11} \end{array}$  \\
\hline
$U_{1,0}$ & $9^2 / 1 \cdot 3$  & $3 \cdot 9 / 1^2$
 & $M_{11}$ \\ \hline \hline
$J'_{2,0}$ & $2^210^2 / 1^25^2$ & $2^210^2 / 1^25^2$
& $J'_{2,0}$ \\
\hline
$\begin{array}{c} L_{1,0} \\ K'_{1,0} \end{array}$ & $2 \cdot 8^2 / 1^24$ & $2
\cdot 8^2 / 1^24$ &
$\begin{array}{c} L_{1,0} \\ K'_{1,0} \end{array}$ \\ \hline
$M_{1,0}$ & $7^2 / 1^2$ & $7^2 / 1^2$ & $M_{1,0}$ \\
\hline
$I_{1,0}$ & $3^26^2 / 1^22^2$ & $3^26^2 / 1^22^2$ & $I_{1,0}$ \\
\hline
\end{tabular}
\end{table}

It turns out that the substitution
$$\mu \mapsto \mu^\flat, \quad c \mapsto c^\flat, \quad H \mapsto H^\flat$$
for the ICIS yields
$$ \mu + \mu^\ast = 24, \quad d = d^\ast, \quad \pi_{X^\ast} = \pi^\ast_X.$$
Moreover, the lattice $H$ admits an embedding into the even
unimodular lattice
$$K_{24}  = (-E_8) \oplus (-E_8) \oplus U \oplus U \oplus U \oplus U$$
of rank $24$.
This lattice can be considered as the full homology lattice of a K3 surface,
$$K_{24} = H_0(S,\ZZ) \oplus H_2(S,\ZZ) \oplus H_4(S,\ZZ),$$
where the inner product on $H_0(S,\ZZ) \oplus H_4(S,\ZZ)$ is defined in such a
way that this lattice corresponds to a unimodular hyperbolic plane $U$. The
orthogonal complement of $H$ is the lattice
$\check{H}$ of the singularity
$X^\ast$ (cf.\ Table~\ref{Table6}).

Let us consider Pinkham's interpretation in the new cases.
The Milnor fibre of a
triangle or quadrilateral isolated hypersurface or complete intersection
singularity can be compactified in such way that after resolving the
singularities one gets a K3 surface $S$ \cite{Pinkham78}. We consider the dual
graph of the curve configuration at infinity in each case.
Let $G(p_1,p_2,p_3,p_4)$ and $\tilde{G}(p_1,p_2,p_3,p_4)$ be the subgraphs of
the graphs of Fig.~\ref{Fig4} and Fig.~\ref{Fig2} respectively obtained by
omitting the vertices $\mu^\flat-1$ and $\mu^\flat$, and $p_1+p_2-1$,
$p_1+p_2+1$, and $\nu$ respectively. Denote by $M(p_1,p_2,p_3,p_4)$ and
$\tilde{M}(p_1,p_2,p_3,p_4)$ the corresponding lattices. Recall that the
homology lattice $H_2(S,\ZZ)$ of the K3 surface is denoted by $L$.

First start with a triangle ICIS $(X,x_0)$ with Dolgachev numbers
$(b_1,b_2,b_3)$. Then the dual graph is the graph $G(b_1,b_2,b_3)$. This yields
an embedding $M(b_1,b_2,b_3) \subset L$ and the orthogonal complement is the
Milnor lattice $H=\tilde{M}(p_1,p_2,p_3,p_4) \oplus U$. By alternative (1) (cf.\
\cite{EW85}) the dual of $(X,x_0)$ is a bimodal series; the Milnor lattice
of the
"virtual" $l=-1$ element of the corresponding series is $M(b_1,b_2,b_3) \oplus
U$. One can even associate a Coxeter element to the dual "virtual"
singularity; it
has order $\mbox{lcm}\,(N_1,N_2)$ where $N_1$, $N_2$ are the degrees of the
equations of $(X,x_0)$ \cite{EW85}, whereas the monodromy operator of $(X,x_0)$
has order $N_2$. There is no Saito duality of characteristic polynomials.
The dual singularities and orders of the monodromy operators are listed in
Table~\ref{Table7b}.

\begin{table}\centering
\caption{The duality: 8 triangle ICIS versus 8 bimodal series}
\label{Table7b}
\begin{tabular}{|c|c|c|c|c|c|c|c|c|}  \hline
Name  & $\mu$  & Dol & Gab & $d$ & $h$ & $h^\ast$ & $\mu^\ast$  & Dual
\\ \hline
$J'_9$ & 9 & 2 3 10 & 2 2 2 3 &  $-4$ & 18 & 144 & 15 & $J_{3,-1}$ \\ \hline
$J'_{10}$ & 10 & 2 4 8 & 2 2 2 4 & 8 & 14 & 84 & 14 & $Z_{1,-1}$ \\ \hline
$J'_{11}$ & 11 & 3 3 7 & 2 2 2 5 & $-12$ & 12 & 60 & 13 & $Q_{2,-1}$ \\ \hline
$K'_{10}$ & 10 & 2 6 6 & 2 3 2 3 & 12 & 12 & 60 & 14  & $W_{1,-1}$ \\ \hline
$K'_{11}$ & 11 & 3 5 5 & 2 3 2 4 & $-20$ & 10 & 40 & 13 & $S_{1,-1}$ \\ \hline
$L_{10}$ & 10 & 2 5 7 & 2 2 3 3 & 11 & 12 & 132 & 14 & $W^\sharp_{1,-1}$ \\
\hline
$L_{11}$ & 11 & 3 4 6 & 2 2 3 4 & $-18$ & 10 & 90 & 13 & $S^\sharp_{1,-1}$ \\
\hline
$M_{11}$ & 11 & 4 4 5 & 2 3 3 3 & $-24$ & 9 & 72 & 13 & $U_{1,-1}$ \\ \hline
\end{tabular}
\end{table}

On the other hand, we can start with a quadrilateral hypersurface singularity
$(X,x_0)$ with Dolgachev numbers $(b_1,b_2,b_3,b_4)$. Then the dual graph is the
graph $G(b_1,b_2,b_3,b_4)$. We obtain an embedding $M(b_1,b_2,b_3,b_4) \subset
L$ and the orthogonal complement is the Milnor lattice of $(X,x_0)$ described
in Table~\ref{Table4}. Here we use alternative (2) for the duality. The reduced
Milnor lattice $H^\flat$ of the dual triangle ICIS according to
Table~\ref{Table4} is the lattice $M(b_1,b_2,b_3,b_4) \oplus U$.

Finally, let $(X,x_0)$ be one of the 5 quadrilateral ICIS. Then the dual graph
is again the graph $G(b_1,b_2,b_3,b_4)$. One has an embedding
$M(b_1,b_2,b_3,b_4) \subset
L$ and the orthogonal complement is the Milnor lattice of $(X,x_0)$ described
in Table~\ref{Table5}. Combining both alternatives (1) and (2), the lattice
$M(b_1,b_2,b_3,b_4)
\oplus U$ can be interpreted as follows: The 5 quadrilateral ICIS are the
initial $l=0$ elements of 8 series of ICIS. To each such series one can again
associate a virtual $l=-1$ element with a well-defined Milnor lattice. Then
$M(b_1,b_2,b_3,b_4) \oplus U$ is the reduced Milnor lattice $\check{H}^\flat$ of
the dual virtual singularity. This correspondence is indicated in
Table~\ref{Table7c}. There is also a duality between the 8 virtual
singularities as indicated in \cite{EW85} (see also
\cite[Table~3.6.2]{Ebeling87}). There is no Saito duality of characteristic
polynomials in both cases. But as we have seen above, using alternative (2) we
get a third correspondence, for which Saito's duality of characteristic
polynomials holds.

\begin{table}\centering
\caption{The duality between the quadrilateral ICIS}
\label{Table7c}
\begin{tabular}{|c|c|c|c|c|c|c|}  \hline
Name & $\mu$ & Dol & Gab & $d$ & $\mu^\ast$ & Dual
\\ \hline
$J'_{2,0}$ & 13 & 2 2 2 6 & 2 2 2 \underline{6} & $-16$ & 11
 & $J'_{2,-1}$ \\ \hline
$L_{1,0}$ & 13 & 2 2 3 5 & $\begin{array}{c} \mbox{2 2 3 \underline{5}} \\
\mbox{2 2 \underline{4} \underline{4}} \end{array}$ & $-28$ &
11
 & $\begin{array}{c} L^\sharp_{1,-1} \\ K^\flat_{1,-1} \end{array}$ \\ \hline
$K'_{1,0}$ & 13 & 2 2 4 4 & $\begin{array}{c} \mbox{2 3 2 \underline{5}} \\
\mbox{2 \underline{4} 2 \underline{4}} \end{array}$ & $-32$ & 11
 & $\begin{array}{c} L_{1,-1} \\ K'_{1,-1} \end{array}$ \\ \hline
$M_{1,0}$ & 13 & 2 3 3 4 &
$\begin{array}{c} \mbox{2 3 \underline{3} \underline{4}} \\
\mbox{2 \underline{3} 3 \underline{4}} \end{array}$ & $-42$ & 11
 & $\begin{array}{c} M^\sharp_{1,-1} \\ M_{1,-1} \end{array}$ \\ \hline
$I_{1,0}$ & 13 & 3 3 3 3 &
$\underline{3}$ $\underline{3}$ $\underline{3}$ $\underline{3}$
 & $-54$ & 11 & $I_{1,-1}$ \\ \hline
\end{tabular}
\end{table}

By Dolgachev's construction \cite{Dolgachev95}, to each case of
Pinkham's construction there corresponds a pair of mirror symmetric families of
K3  surfaces. Moreover, also to each case where we only have a pair of
lattices embedded as orthogonal complements to each other in the lattice
$K_{24}$ (cf.\ Table~\ref{Table6}) there corresponds such a mirror pair.

One can also investigate Kobayashi's duality of weight systems for our
extension of Arnold's strange duality. As already observed by Kobayashi
\cite{Kobayashi95}, only some of the weight systems of the quadrilateral
hypersurface singularities have dual weight systems, the dual weight systems
are in general not unique and they correspond again to isolated
hypersurface singularities. Since an ICIS has two degrees $N_1$ and $N_2$, it is
not quite clear how to generalize the notion of a weighted magic square. One
possibility would be to work with the sum of the degrees $N := N_1 + N_2$ and
to use also $3 \times 4$ and $4 \times 4$ matrices instead of $3 \times 3$
matrices. Then one finds again that in some cases there does not exist a dual
weight system, the dual weight systems are in general not unique, most cases are
self-dual, and only in the cases
$J'_{10}
\leftrightarrow Z_{1,0}$,
$K'_{11} \leftrightarrow S_{1,0}$, $J'_{2,0} \leftrightarrow J'_{2,0}$, and
$M_{1,0} \leftrightarrow M_{1,0}$ of our duality there exist weighted magic
squares giving a duality of the corresponding weight systems.

However, there is a relation between our extended duality and a polar duality
between the Newton polytopes generalizing Kobayashi's observation for Arnold's
strange duality. This can be used to explain Saito's duality of characteristic
polynomials. For details see the forthcoming paper \cite{Ebeling98}.

\section{Singular moonshine}

Let us consider the symbols $\pi\pi^\ast$ of Table~\ref{Table7}. It turns out
that they all occur in the list of Kondo, too. These pairs and the pairs from
the original Arnold duality correspond to self-dual Frame shapes of the group
$G$ with trace $-2$, $-3$, or
$-4$. By examining Kondo's list one finds that there are 22 such Frame shapes
and all but 3 occur. They are listed in Table~\ref{Table8}. Here we use the
ATLAS notation \cite{ATLAS} for the conjugacy classes. For each value of the
trace one symbol is missing.

\begin{table}\centering
\caption{Self-dual Frame shapes of $G$ with trace $-2$, $-3$ or $-4$}
\label{Table8}
\begin{tabular}{|c|c|c|c|c|}  \hline
ATLAS & Frame  &  & Niemeier & Duality
\\ \hline
21A & $2^23^27^242^2 / 1^26^214^221^2$ &  $*$ &
 & $E_{12} \leftrightarrow E_{12}$
\\ \hline
15E & $2^23 \cdot5 \cdot 30^2 / 1^26 \cdot 10 \cdot 15^2$  & $*$ &
$D_{16} \oplus E_8$
 & $E_{13} \leftrightarrow Z_{11}$
\\ \hline
24B & $2 \cdot 3^24 \cdot 24^2 / 1^26 \cdot 8^2 12$ & $*$ &
 & $E_{14} \leftrightarrow Q_{10}$
\\ \hline
11A & $2^222^2 / 1^211^2$  & $*$ & $D_{12}^2$
 & $Z_{12} \leftrightarrow Z_{12}$
\\ \hline
18B & $2 \cdot 3 \cdot 18^2 / 1^26 \cdot 9$  & $*$ &
$A_{17} \oplus E_7$
 & $Z_{13} \leftrightarrow Q_{11}$
\\ \hline
15B & $3^215^2 / 1^25^2$  & $*$ &
 & $Q_{12} \leftrightarrow Q_{12}$
\\ \hline
20A & $2^25^220^2 / 1^24^210^2$  & $*$ &
 & $W_{12} \leftrightarrow W_{12}$
\\ \hline
16B & $2 \cdot 16^2 / 1^28$  & $*$ & $A_{15} \oplus D_9$
 & $W_{13} \leftrightarrow S_{11}$
\\ \hline
13A & $13^2 / 1^2$  & $*$ & $A^2_{12}$
 & $S_{12} \leftrightarrow S_{12}$
\\ \hline
12E & $4^212^2 / 1^23^2$  & $*$ &
 & $U_{12} \leftrightarrow U_{12}$
\\ \hline
9C & $2^33^218^3 /1^36^29^3$ &  & $D_{10} \oplus E_7^2$
 & $J'_9 \leftrightarrow J_{3,0}$
\\ \hline
7B & $2^314^3 / 1^37^3$  & $*$ & $D^3_8$
 & $J'_{10} \leftrightarrow Z_{1,0}$
\\ \hline
12K & $2^23 \cdot 12^3 / 1^34 \cdot 6^2$  & $*$ & $A_{11}
\oplus D_7
\oplus E_6$
 & $\begin{array}{c} J'_{11} \leftrightarrow Q_{2,0} \\
                     K'_{10} \leftrightarrow W_{1,0} \\
                     L_{10} \leftrightarrow W_{1,0}
\end{array}$
\\ \hline
10E & $2 \cdot 10^3 / 1^35$  & $*$  & $A_9^2
\oplus D_6$
 & $\begin{array}{c} K'_{11} \leftrightarrow S_{1,0} \\
                     L_{11} \leftrightarrow S_{1,0}
\end{array}$
\\ \hline
9A & $9^3 / 1^3$  & & $A_8^3$
 & $M_{11} \leftrightarrow U_{1,0}$
\\ \hline
5B & $2^410^4 / 1^45^4$
&  & $D_6^4$
 & $J'_{2,0} \leftrightarrow J'_{2,0}$
\\ \hline
8C & $2^28^4 / 1^44^2$  & & $A_7^2 \oplus D_5^2$
 & $\begin{array}{c} L_{1,0} \leftrightarrow L_{1,0} \\
                     L_{1,0} \leftrightarrow K'_{1,0} \\
                     K'_{1,0} \leftrightarrow K'_{1,0}
\end{array}$
\\ \hline
7A & $7^4 / 1^4$  & & $A_6^4$
 & $M_{1,0} \leftrightarrow M_{1,0}$
\\ \hline
6A & $3^46^4 / 1^42^4$  & &
 & $I_{1,0} \leftrightarrow I_{1,0}$
\\ \hline
10A & $5^210^2 / 1^22^2$  &
$*$ & &
\\ \hline
15A & $2^33^35^330^3 / 1^36^310^315^3$  &
$*$ & $E_8^3$ &
\\ \hline
12A & $2^43^412^4 /1^44^46^4$  &  & $E_6^4$ &
\\ \hline
\end{tabular}
\end{table}

Special elements of $G$ correspond to the deep holes of the Leech lattice. A
{\em deep hole} of the Leech lattice $\Lambda_{24}$ is a point in $\RR^{24}$
which has maximal distance from the lattice points. It is a beautiful theorem of
J.~H.~Conway, R.~A.~Parker, and N.~J.~A.~Sloane (\cite{CPS82}, see also
\cite{Ebeling94}) that there are 23 types of deep holes in $\Lambda_{24}$ which
are in one-to-one correspondence with the 23 isomorphism classes of even
unimodular lattices in
$\RR^{24}$ containing roots, which were classified by H.-V.~Niemeier
\cite{Niemeier73}. These lattices are characterized by the root systems which
they contain. The Frame shape of the Coxeter element of such a root system is
also the Frame shape of an automorphism of the Leech lattice. We have
indicated in Table~\ref{Table8} the type of the root system of the Niemeier
lattice, if a Frame shape corresponds to the Coxeter element of such a root
system.

The automorphism group $G$ of the Leech lattice contains the Mathieu group
$M_{24}$ (see e.g.\ \cite{Ebeling94}). S.~Mukai \cite{Mukai88} has
classified the
finite automorphism groups of K3 surfaces (automorphisms which leave the
symplectic form invariant) and shown that they admit a certain embedding into
the Mathieu group $M_{24}$. He gives a list of 11 maximal groups such that
every finite automorphism group imbeds into one of these groups. A table of the
centralizers of the conjugacy classes of $G$ can be found in \cite{Wilson83}.
For an element $g
\in G$, denote its centralizer by $Z(g)$ and the finite cyclic group
generated by
$g$ by
$\langle g \rangle$. We have marked by ($\ast$) in Table~\ref{Table8} the cases
where there is an obvious inclusion of
$Z(g) / \langle g \rangle$ in one of Mukai's groups. It follows that in these
cases there is a K3 surface with an operation of $Z(g) / \langle g \rangle$ on
it by symplectic automorphisms.

To a Frame shape
$$\pi = \prod_{m | N} m^{\chi_m}$$
one can associate a modular function \cite{Kondo85}. Let
$$\eta(\tau) = q^{1/24} \prod^\infty_{n=1} (1-q^n), \quad q=e^{2\pi i \tau},
\tau \in \HH,$$
be the Dedekind $\eta$-function. Then define
$$\eta_\pi(\tau) = \prod_{m / N} \eta(m\tau)^{\chi_m}.$$
Saito \cite{Saito94} has proved the identity
$$\eta_\pi \left( - \frac{1}{N\tau} \right) \eta_{\pi^\ast}(\tau) \sqrt{d} =
1,$$
where $d = \prod m^{\chi_m}$ and $\pi^\ast$ is the dual Frame shape. From this
it follows that $\eta_{\pi\pi^\ast}$ is a modular function for the group
$$\Gamma_0(N) = \left\{ \left( \begin{array}{cc} a & b \\ c & d \end{array}
\right) \in \mbox{SL}_2(\ZZ) \left| \right. c \equiv 0 (N) \right\}.$$

\addvspace{3mm}

\noindent {\bf Question~1} Let $\pi\pi^\ast$ be one of the self-dual Frame
shapes of Table~\ref{Table8}. Is there any relation of $\eta_{\pi\pi^\ast}$ to
the analogue of Dedekind's eta function for K3 surfaces considered in
\cite{JT95} ?

\addvspace{3mm}

The Frame shape $\pi\pi^\ast$ is the Frame shape of the
operator $c \oplus c^\ast$ which can be considered as an automorphism of a
sublattice of finite index of the even unimodular 24-dimensional lattice
$K_{24}$, which is the full homology lattice of a K3 surface. The
lattice $K_{24}$ has the same rank as the Leech lattice, but contrary to the
Leech lattice it is indefinite and has signature
$(4,20)$.

\addvspace{3mm}

\noindent {\bf Question~2 } Is there an explanation for this strange
correspondence between operators of different lattices?

\addvspace{3mm}

Is it only a purely combinatorial coincidence? One can try to classify finite
sequences $(\chi_1, \chi_2, \ldots , \chi_N)$ with the following properties:
\begin{itemize}
\item[(1)] $\chi_m \in \ZZ$ for all $m=1, \ldots , N$,
\item[(2)] $\sum m\chi_m = 24$,
\item[(3)] $\chi_m = 0$ for $m \! \! \not| N$,
\item[(4)] $\chi_m = - \chi_{N/m}$ for $m | N$,
\item[(5)] $\prod m^{\chi_m} \in \NN$,
\item[(6)] $\chi_1  \in \{-2,-3,-4\}$,
\item[(7)] $|\chi_m| \leq |\chi_N|$ for $m | N$.
\end{itemize}
By a computer search, one finds for $N \leq 119$ in addition to the 22 Frame
shapes of Table~\ref{Table8} only the following Frame shapes:
$$3^26 \cdot 12^2 /1^22 \cdot 4^2, \quad 2^43^44^424^4 / 1^46^48^412^4, \quad
2^24^25^240^2 / 1^28^210^220^2.$$
These Frame shapes also appear in Kondo's tables, namely in the tables of
certain transforms of the Frame shapes of $G$ \cite[Table~III, 4C; Table~II,
12A; Table~II, 20A]{Kondo85}.

\addvspace{10mm}

{\sc Institut f\"{u}r Mathematik, Universit\"{a}t Hannover, Postfach 6009,
D-30060 Hannover, Germany}

{\em E-mail address}: ebeling@math.uni-hannover.de
\end{document}